\documentclass{amsart}
\usepackage[utf8]{inputenc}
\usepackage{amsmath,amsthm,amsfonts,amssymb,mathtools}
\usepackage{mathrsfs}
\makeatletter
\@namedef{subjclassname@2020}{\textup{2020} Mathematics Subject Classification}
\makeatother

\usepackage{hyperref}
\usepackage{capt-of}
\usepackage{relsize}

\newcommand{\Be}{\begin{equation}}
\newcommand{\Ee}{\end{equation}}
\newcommand{\Bea}{\begin{eqnarray}}
\newcommand{\Eea}{\end{eqnarray}}
\newcommand{\Bel}{\begin{align}}
\newcommand{\Eel}{\end{align}}
\newcommand{\Beas}{\begin{eqnarray*}}
\newcommand{\Eeas}{\end{eqnarray*}}
\newcommand{\Benu}{\begin{enumerate}}
\newcommand{\Eenu}{\end{enumerate}}
\newcommand{\Bi}{\begin{itemize}}
\newcommand{\Ei}{\end{itemize}}
\newcommand\supp{\operatorname{supp}}

\def\R{{\mathbb R}}
\def\Z{{\mathbb Z}}
\def\N{{\mathbb N}}

\theoremstyle{plain}
\newtheorem{thm}{Theorem}[section]

\newtheorem{lem}[thm]{Lemma}
\newtheorem{prop}[thm]{Proposition}

\newtheorem{conj}{Conjecture}
\theoremstyle{remark}
  
\theoremstyle{definition}

\numberwithin{equation}{section}

\newcommand{\RNum}[1]{\uppercase\expandafter{\romannumeral #1\relax}}

\newcommand{\cA}{\mathcal A}

\usepackage{wrapfig}
\usepackage{tikz}

\newcommand{\vl}{{\vec\ell}}

\makeatletter
\@namedef{subjclassname@2020}{\textup{2020} Mathematics Subject Classification}
\makeatother

\title{Maximal estimates for averages over degenerate hypersurfaces}

\author[Sewook Oh]{Sewook Oh}

\address{June E Huh Center for Mathematical Challenges, Korea Institute for Advanced Study, Seoul 02455, Republic of Korea}
\email{sewookoh@kias.re.kr}

\begin{document}

\maketitle

\begin{abstract}
We study $L^p$ boundedness of the maximal average over dilations of a smooth hypersurface $S$. When the decay rate of the Fourier transform of a measure on $S$ is $1/2$, we establish the optimal maximal bound, which settles the conjecture raised by Stein. Additionally, when $S$ is not flat, we verify that the maximal average is bounded on $L^p$ for some finite $p$, which generalizes the result by Sogge and Stein.
\end{abstract}

\section{Introduction}

Let $S$ be a smooth hypersurface in $\R^d$ which is given by the graph 
 of a smooth function $\gamma:\R^{d-1}\rightarrow\R$.  We define a measure $\sigma_t$ supported on $tS$ by 
\[
\langle \sigma_t,f\rangle =\int f(t\Gamma(u))\psi(u) du,
\]
where $\Gamma(u)=(u,\gamma(u))$ and $\psi\in C^\infty_0(\R^{d-1})$ is a nonnegative cut-off function whose support is in  $B(0,1):=\{u\in\R^{d-1}:|u|\le1\}$.
We then define an averaging operator $\mathcal A$ as
\[
\mathcal Af(t,x)=f\ast \sigma_t(x).
\]
In this paper, we are concerned with the $L^p$ boundedness of the maximal operator $\mathcal M$ associated to a smooth hypersurface $S$, which is given by
\[
\mathcal Mf(x)=\sup_{t>0}|\mathcal Af(t,x)|.
\]

The $L^p$ bounds of maximal operators associated to hypersurfaces have been extensively studied since the 1970s, generating an immense body of literature (see \cite{S, IKM2} and references therein). The study of maximal averages has been initiated by the celebrated Stein's spherical maximal theorem \cite{S2}, which tells that $\mathcal M$ is bounded on $L^p$ for $p>d/(d-1)$ when $S$ is a sphere in $\R^d$ for $d\ge 3$. Later, Bourgain \cite{B} extended this result to the case $d=2$. In \cite{S2}, the decay rate of the Fourier transform of $\sigma_1$ plays a crucial role. The connection between Fourier decay rate and $L^p$ maximal bounds has been explored in various studies (\cite{G,R,CM}). By Greenleaf \cite{G} and Rubio de Francia \cite{R}, it is known that the maximal operator $\mathcal M$ is bounded on $L^p$ for $p>1+1/(2q)$, if  
\Be
\label{Fourierdecay}
|\hat\sigma_1(\xi)|\le C|\xi|^{-q},
\Ee
for $q>1/2$ and some $C>0$. This result provides the optimal $L^p$ maximal bounds when $S$ has nonvanishing Gaussian curvature, specifically for $q=(d-1)/2$. However, it is not optimal for $1/2<q<(d-1)/2$. Notably, as discussed in the literature (see \cite{BIM,IS1,IS2,ISS,NSW}), Fourier decay alone is not sufficient to characterize the $L^p$ boundedness of the maximal operator when $q>1/2$.

In contrast, for $q\le 1/2$, Fourier decay becomes a more significant factor in determining the $L^p$ boundedness of the maximal operator.
Stein, for $q=1/2$, Iosevich--Sawyer, for $0<q<1/2$, conjectured the following.
\begin{conj}
Assume that \eqref{Fourierdecay} holds for $0<q\le 1/2$. Then $\mathcal M$ is bounded on $L^p$ for $p>1/q$.
\end{conj}
For the case $q\le 1/2$, Stein's $L^2$- method no longer works, which makes this problem poorly understood. The conjecture is currently known  to hold only for specific cases: 
 When $S$ has at least one non-vanishing principal curvature everywhere, as shown by Sogge \cite{Sogge}, for $d=2$, as proved by Iosevich \cite{I}, and in the case $d=3$, as established by Ikromov--Kempe--M\"uller \cite{IKM1,IKM2} (see \cite{KLO2,BGHS,KLO} for results when $S$ is a curve). In these studies, it is clear that the conjectured range, $p>1/q$, is optimal.
 However, when $d\ge4$, the optimal result for general hypersurfaces remains unknown.
 A partial result was known by Sogge and Stein \cite{SS1}, which verifies nontrivial $L^p$ bounds for maximal averages when  $p$ is sufficiently large, provided that the Gaussian curvature of $S$ does not vanish to infinite order at any point of  $S$.

In this paper,  we assume that $\gamma$ is of finite type, meaning that for all $u\in B(0,1)$, there exists a multi-index $\alpha=\alpha(u)\in \N_0^{d-1}$ such that $\partial^\alpha\gamma(u)\neq0$ and $|\alpha|\ge2$. Here, $\N_0:=\N\cup \{0\}$ and $|\alpha|=\sum_i|\alpha_i|$. Indeed,
if $\gamma$ is flat at some point $u_0$ satisfying $\psi(u_0)\neq0$, i.e., $\partial^\alpha\gamma(u_0)=0$ for all  $|\alpha|\ge2$,  then \eqref{Fourierdecay} fails for all $q>0$.\footnote{\eqref{Fourierdecay} with $q>0$ implies that $\gamma$ is not flat at $u_0$. Indeed, as discussed in \cite{CCW},  \eqref{Fourierdecay} with $0<q<1$ implies the sublevel set estimate, $\int_{\{u:|\gamma(u)-\gamma(u_0)-\nabla\gamma(u_0)\cdot(u-u_0)|\le \delta\}}\psi(u)du\le C\delta^q$. By Taylor's theorem, it implies that $\partial^\alpha \gamma(u_0)\neq0$ for some $|\alpha|\ge2$.}  Thus, this case is of no interest.

Under the finite type assumption, we prove the optimal $L^p$ maximal bound when $q=1/2$, which confirms Stein's conjecture.
\begin{thm}\label{thm:1/2}
    Let $\gamma:\R^{d-1}\rightarrow\R$ be a smooth function of finite type. Suppose that \eqref{Fourierdecay} holds for $q=1/2$. Then $\mathcal M$ is bounded on $L^p$ for $p>2$.
\end{thm}
When $0<q<1/2$, we extend the result by Sogge and Stein \cite{SS1} to general hypersurface, removing the assumption they required.
\begin{thm}\label{thm:<1/2}
    Let $\gamma:\R^{d-1}\rightarrow\R$ be a smooth function of finite type. Suppose that \eqref{Fourierdecay} holds for $0<q<1/2$. Then there exists $p_0>2$ such that $\mathcal M$ is bounded on $L^p$ for $p>p_0$.
\end{thm}

To prove two theorems, we use local smoothing estimates. The local smoothing estimate for averages is given by the following inequality:
\Be
\label{localsmoothing}
    \|\chi \mathcal Af\|_{L^p_{\mathfrak b}(\R^{d+1})}\le C \|f\|_{L^p(\R^d)},
\Ee
where $\chi\mathcal A f(t,x)=\chi(t)\mathcal Af(t,x)$ and $\chi$ is a smooth function supported on $(1/2,4)$.
This estimate has  direct applications to $L^p$ maximal bounds (see \cite{MSS,SS2,PS,BGHS,KLO}). It is well-known that if \eqref{localsmoothing} holds with $\mathfrak b>1/p$, then the maximal operator $\mathcal M$ is bounded on $L^p$. Therefore, the above theorems are direct consequences of the following local smoothing estimate.
\begin{thm}\label{thm:main}
    Let $k\ge2$ and $\gamma:\R^{d-1}\rightarrow\R$ be a smooth function. Suppose that there exists a multi-index $\alpha$ with $|\alpha|=k$ such that  $\partial^\alpha\gamma(u)\neq 0$ for  all $u\in B(0,1)$. Then, for $p\ge \max\{4k-4,6\}$ and $\epsilon>0$,
    \eqref{localsmoothing} holds with $\mathfrak b=2/p-\epsilon$.
\end{thm}
After a suitable decomposition,
by smoothness of $\gamma$ and compactness of  the support of $\psi$, Theorem \ref{thm:<1/2} follows directly from Theorem \ref{thm:main}. When $q=1/2$, \eqref{Fourierdecay} provides the $L^2$ estimate, \eqref{localsmoothing} for $p=2,\mathfrak b=1/2$. By interpolation between this $L^2$ estimate and Theorem \ref{thm:main}, we  observe that \eqref{localsmoothing} holds for some $\mathfrak b>1/p$ when $p>2$, thus confirming Theorem \ref{thm:1/2}. 

The smoothing order $2/p$ in Theorem \ref{thm:main} is sharp in the sense that there exists a smooth function $\gamma$ satisfying the conditions in Theorem \ref{thm:main} such that \eqref{localsmoothing} fails for $\mathfrak b>2/p$. For example, when $\gamma(u_1,\cdots, u_{d-1})=1+u_1^k$, a modification of the counter example in \cite{KLO} tells that \eqref{localsmoothing} is valid only if $\mathfrak b\le 2/p$.

In previous studies (\cite{W, PS,KLO}), decoupling inequalities play a crucial role in establishing local smoothing estimates for large $p$. These works reveal that averaging operators associated to nondegenerate curves are closely related to the decoupling inequality for the conic extension of the moment curve. However, when $S$ is a degenerate hypersurface, the connection between averaging operators and the corresponding decoupling inequality is not well understood, primarily due to the difficulty of handling the degeneracy of $S$.

For curves, a common strategy for addressing degeneracy involves decomposition away from degeneracy combined with scaling. For hypersurfaces in higher dimensions, however, identifying suitable scaling becomes significantly more challenging. To overcome these difficulties, we focus on analyzing  a set away from degeneracy, defined as
\[\{u\in \R^{d-1}:\rho\le(\sum_{|\alpha|=2}|\partial^\alpha\gamma(u)|^2)^{1/2}  \le4\rho \}.\]
We show that this set can be covered by a collection of balls on which the higher order terms of $\gamma$ can be treated as error terms. The number of balls in this collection and their respective volumes are carefully controlled (see Proposition \ref{prop:il}), allowing us to reduce the analysis to small, localized pieces of $\Gamma$.

 For each small component, we further observe that, after suitable decomposition and projection, the essential Fourier support of the average over a piece of $\Gamma$ is contained in a neighborhood of a cone in $\R^3$  (see Proposition \ref{prop:projcone}). This geometric observation enables us to apply the well-established decoupling inequality for cones, thereby achieving the desired local smoothing estimate.

\subsection*{Organization of the paper}
In section \ref{sec:red}, the proof of Theorem \ref{thm:main} is reduced to establishing Proposition \ref{prop:a1}. Proposition \ref{prop:a1} is proved in section \ref{sec:decomp}, while assuming that Theorem \ref{thm:avl} holds. The proof of Theorem \ref{thm:avl} is provided in section \ref{sec:decoup}.


\section{Frequency localization}\label{sec:red}

Let $\mathcal B\ge1$ be a large number and $N$ be the smallest integer which is greater than $(2C_\ast)/\epsilon$. The constant $C_\ast$ will be chosen later, but note that $C_\ast$ depends only on $d$.
From now on, we assume that  
\Be
\label{cond:gamma}
 \sum_{\alpha\in \N_0^{d-1}:|\alpha|= k}|\partial^{\alpha}\gamma(u)|\ge \mathcal B^{-1} \text{ and } \sup_{|\alpha|\le N+1}|\partial^\alpha\gamma(u)|\le \mathcal B,
\Ee
for  $u\in B(0,1)$. 

\subsection*{Notation}
For nonnegative quantities $A$ and $B$, we denote $A\lesssim B$ if there exists a constant $C$ depending only on $d,p,k,\epsilon, \mathcal B$ such that $A\le CB$, and $A\sim B$ means that $A\lesssim B$ and $B\lesssim A$. By $A=O(B)$, we denote $|A|\lesssim B$.

We define an averaging operator with a symbol $a\in C([1/2,4]\times \R^{d-1}\times \R\times \R^{d})$ by
\[
\mathcal A[a]f(t,x):=(2\pi)^{-d-1}\iint\big(\iint e^{-it'(\tau+\Gamma(u)\cdot \xi)} a(t',u,\tau,\xi)dudt'\big)e^{i(x\cdot\xi+t\tau)}\hat f(\xi)d\xi d\tau.
\]
Note that  $\mathcal A[\chi\psi]f=\chi\mathcal Af$. 
Set $\lambda \ge1$ and
\[
\mathbb A_\lambda=\{(\tau,\xi)\in \R\times \R^d:\lambda/2\le|(\tau,\xi)|\le 2\lambda\}.
\]
Theorem \ref{thm:main} is a consequence of the following theorem.
\begin{thm}\label{thm:main'}
Let $\lambda\ge 1$.
    Suppose that a symbol $a$ satisfies $\supp a\in [1/2,4]\times B(0,1)\times \mathbb A_\lambda$ and
    \Be
    \label{cond:a}
    |\partial_t^l\partial_u^\alpha\partial_{(\tau,\xi)}^\beta a(t,u,\tau,\xi)|\lesssim |(\tau,\xi)|^{-|\beta|}
    \Ee
    for $(l,\alpha,\beta)\in \N_0\times \N_0^{d-1}\times\N_0^{d+1}$ satisfying $l,|\alpha|\le N,$ and $|\beta|\le d+2$. Then for $p>\max\{4k-4,6\}$ and $\epsilon>0$ such that 
    \Be
    \label{ineq:main}
    \|\mathcal A[a]f\|_{L^p(\R^{d+1})}\lesssim \lambda^{-(\frac 2p-\epsilon)}\|f\|_p.
    \Ee
\end{thm}


From now on, we fix a small constant $\epsilon_1:=\epsilon/C_\ast$. 
To show Theorem \ref{thm:main'}, we decompose a symbol $a$ into two cases. The first case is either:
\Be
\label{easycase}
|(\tau,\xi)\cdot(1,\Gamma(u))|\ge \lambda^{\epsilon_1} \text{ or } |\xi\cdot\partial_j\Gamma(u)|\ge \lambda^{\frac 12+\epsilon_1}
\Ee
for some $j=1,\cdots,d-1$. Here, by $\partial_j$ we denote $\partial_{u_j}$. The second case is the remaining case:
\Be
\label{maincaes}
|(\tau,\xi)\cdot(1,\Gamma(u))|\le 2\lambda^{\epsilon_1} \text{ and } |\xi\cdot\partial_j\Gamma(u)|\le 2\lambda^{\frac 12+\epsilon_1}
\Ee
for all $j=1,\cdots, d-1$.
For decomposing a symbol $a$, define a smooth function $\chi_\lambda:\R^{d-1}\times \R\times \R^d\rightarrow \R$ by
\[
\chi_\lambda(u,\tau,\xi)=\eta_0(\lambda^{-\epsilon_1}(\tau,\xi)\cdot(1,\Gamma(u)))\prod_{j=1}^{d-1}\eta_0(\lambda^{-1/2-\epsilon_1}\xi\cdot\partial_j\Gamma(u)),
\]
where $\eta_0:\R\rightarrow \R$ is a smooth function satisfying that $\eta_0\equiv 1$ on $[-1,1]$, $\supp \eta_0\subset[-2,2]$.
Then we can easily check that \eqref{easycase} holds on $\supp (a(1-\chi_\lambda))$ and \eqref{maincaes} holds on $\supp (a\chi_\lambda)$.
By triangle inequality, proving Theorem \ref{thm:main'} reduces to verifying the following proposition. 
\begin{prop}\label{prop:a1}
    Let $\epsilon_1>0$. Suppose that a symbol $a$ satisfies the same assumption in Theorem \ref{thm:main'}. Then we have 
    \Be
    \label{aa01}
    \|\mathcal A[a(1-\chi_\lambda)]f\|_{L^p(\R^{d+1})}\lesssim \lambda^{-\frac{2}p}\|f\|_p,
    \Ee
    for $2\le p\le\infty$ and 
    \Be
    \label{aa00}
    \|\mathcal A[a\chi_\lambda]f\|_{L^p(\R^{d+1})}\lesssim \lambda^{-\frac{2}p+C_\ast\epsilon_1}\|f\|_p,
    \Ee
    for $p>\max\{4k-4,6\}$.
\end{prop}

In the rest of this section, we consider the easier case $\mathcal A[a(1-\chi_\lambda)]$. Indeed,
since we have \eqref{easycase}  for this case, we can easily prove \eqref{aa01} using integration by parts. The proof of \eqref{aa00} is provided in Section \ref{sec:decomp}.

\subsection{Proof of \eqref{aa01}}
In this subsection, for simplifying the notation, by $a_\lambda$ we denote $a(1-\chi_\lambda)$. 
To show \eqref{aa01} for $2\le p\le \infty$, we investigate $\|\cA[a_\lambda]\|_{L^p(\R^{d+1})}$ for $p=2,\infty$.  By interpolation between $p=2$ and $p=\infty$, it is sufficient to prove \eqref{aa01}
    for $p=2,\infty$.

    \subsection*{When p=2}
    Plancherel's theorem tells that
    \[
    \|\cA[a_\lambda]f\|_{L^2(\R^{d+1})}^2\sim \iint \big(\iint e^{-it'(\tau+\Gamma(u)\cdot \xi)} a_\lambda(t',u,\tau,\xi)dudt'\big)^2|\hat f(\xi)|^2d\xi d\tau.
    \]
    Since $\supp a\subset [1/2,4]\times B(0,1)\times \mathbb A_\lambda$, we can observe that    
    $\supp_\tau a_\lambda\subset [-2\lambda,2\lambda]$. From this observation and Plancherel's theorem, \eqref{aa01} for $p=2$ is reduced to prove a suitable bound on the Fourier multiplier: for $n=0,1$
    \Be
    \label{aa1p2}
    \big|\iint e^{-it'(\tau+\Gamma(u)\cdot \xi)} a_{\lambda,n}(t',u,\tau,\xi)dudt'\big|\lesssim \lambda^{-\frac 32},
    \Ee
    where 
    $a_\lambda=a_{\lambda,0}+a_{\lambda,1}$ and 
    \[
    a_{\lambda,0}(t,u,\tau,\xi)=a_\lambda(t,u,\tau,\xi)\big(1-\eta_0(4\lambda^{-\epsilon_1}(\tau,\xi)\cdot(1,\Gamma(u)))\big).
    \]
    By definition, one can check that $|(\tau,\xi)\cdot(1,\Gamma(u))|\gtrsim \lambda^{\epsilon_1}$ on $\supp a_{\lambda,0}$. Using this lower bound on phase, we can easily verify \eqref{aa1p2} for $n=0$. Indeed, after applying integration by parts $N$-times with respect to $t'$, we have that the left-hand side of \eqref{aa1p2} for $n=0$ is bounded by 
    \[
    \big|\iint e^{-it'(\tau+\Gamma(u)\cdot \xi)} \frac{\partial_{t'}^Na_{\lambda,0}(t',u,\tau,\xi)}{((\tau,\xi)\cdot(1,\Gamma(u)))^N}dudt'\big|\lesssim \lambda^{-\frac 32} \sup_{t',u}|\partial_{t'}^N a_{\lambda,0}|,
    \]
    since $\supp_{t',u} a_{\lambda,0}\subset [1/2,4]\times B(0,1)$ and $N\epsilon_1>3/2$. Here, we denote $\supp_{t',u}a=\cup_{\tau,\xi}\supp a(\cdot,\cdot,\tau,\xi)$. 
    We also use the same convention with other variables. To prove \eqref{aa1p2} for $n=0$, it is enough to show that   
     $|\partial_{t}^N a_{\lambda,0}|\lesssim 1$. By definitions of $a_{\lambda}$ and $a_{\lambda,0}$, we can check that $|\partial_{t'}^N a_{\lambda,0}|\le|\partial_{t'}^N a|$. Thus, the assumption \eqref{cond:a} implies the desired.

    Now, we focus on \eqref{aa1p2} for $n=1$.
    Since \eqref{easycase} holds on $\supp a_\lambda$ and $|(\tau,\xi)\cdot(1,\Gamma(u))|\le \lambda^{\epsilon_1}/2$ on $\supp a_{\lambda,1}$, we can observe that
    \Be
    \label{cond:ae1}
    |\nabla_u(\xi\cdot \Gamma(u))|\gtrsim \lambda^{1/2+\epsilon_1}
    \Ee
    on $\supp a_{\lambda,1}$.
    For verifying \eqref{aa1p2} with $n=1$, we use integration by parts with respect to
        \[
    D=\frac{\nabla_u(\xi\cdot \Gamma)\cdot \nabla_u}{|\nabla_u(\xi\cdot \Gamma)|^2}.
    \]
    Since $e^{-it'(\tau+\Gamma(u)\cdot \xi)}=(-it')^{-N}D^N(e^{-it'(\tau+\Gamma(u)\cdot \xi)})$, applying integration by parts $N$-times yields that the left-hand side of \eqref{aa1p2} for $n=1$ is bounded by 
    \[
    \big|\iint e^{-it'(\tau+\Gamma(u)\cdot \xi)} (D^\intercal)^Na_{\lambda,1}(t',u,\tau,\xi)dudt'\big|\lesssim  \sup_{t',u}|(D^\intercal)^N a_{\lambda,1}|.
    \]
    Here, we also use $\supp_{t',u} a_{\lambda,1} \subset [1/2,4]\times B(0,1) $. Thus, to establish \eqref{aa1p2} for $n=1$, it is sufficient show that
     \Be
    \label{upbd:ae1}
    |(D^\intercal)^N a_{\lambda,1}|\lesssim \lambda^{-\epsilon_1N}\le \lambda^{-3/2}. 
    \Ee
    
    For  \eqref{upbd:ae1},
    we first check suitable estimates for derivatives of $a_{\lambda,1}$. Recall the definition of $a_{\lambda,1}$,
    \[
    a_{\lambda,1}=a(1-\chi_\lambda)\eta_0(4\lambda^{-\epsilon_1}(\tau,\xi)\cdot(1,\Gamma(u))).
    \]
    To estimate the derivatives of $a_{\lambda,1}$, we use upper bounds for the derivatives of $(\tau,\xi)\cdot(1,\Gamma(u))$. Considering $\supp_{\tau,\xi} a_{\lambda,1}\subset  \mathbb A_\lambda$,  we can observe that
    \Be
    \label{ineq:gamma1}
    \lambda^{-\epsilon_1}|(\tau,\xi)\cdot\partial_u^\alpha(1,\Gamma(u))|\lesssim \big(\lambda^{-\epsilon_1}|\nabla_u(\xi\cdot \Gamma(u))|\big)^{|\alpha|},
    \Ee
on $\supp a_{\lambda,1}$, for $|\alpha|\le N+1$ and
    \Be
    \label{ineq:gamma2}
    \lambda^{-\frac 12-\epsilon_1}|\xi\cdot\partial_u^\alpha\Gamma(u)|\lesssim \big(\lambda^{-2\epsilon_1}|\nabla_u(\xi\cdot \Gamma(u))|\big)^{|\alpha|-1},
    \Ee
on $\supp a_{\lambda,1}$, for $2\le|\alpha|\le N+1$. Here, we utilize the fact that $\lambda^{1/2}\lesssim \lambda^{-\epsilon_1}|\nabla_u(\xi\cdot\Gamma(u))|$ as stated in \eqref{cond:ae1}. Using \eqref{ineq:gamma1} and \eqref{ineq:gamma2},  a computation gives that    
\[
|\partial_u^\alpha\big( (1-\chi_\lambda)\eta_0(4\lambda^{-\epsilon_1}(\tau,\xi)\cdot(1,\Gamma(u)))\big)|\lesssim \big(\lambda^{-\epsilon_1}|\nabla_u(\xi\cdot\Gamma(u))|\big)^{|\alpha|}
\]
    on $\supp a_{\lambda,1}$, for $|\alpha|\le N$. 
    Now, the above estimate combined with 
     \eqref{cond:a}
    yields that for $|\alpha|\le N$,
    \Be
    \label{ineq:paa1}
    |\partial_u^\alpha a_{\lambda,1}|\lesssim \big(\lambda^{-\epsilon_1}|\nabla_u(\xi\cdot\Gamma(u))|\big)^{|\alpha|}.
    \Ee
    Additionally, from \eqref{cond:ae1} and \eqref{ineq:gamma2}, we can see that for $|\alpha|\le N$, $j=1,\cdots, d-1$,
    \Be
    \label{ineq:pgamma}
    \big|\partial_u^\alpha(\frac{\partial_j(\xi\cdot\Gamma)}{|\nabla_u(\xi\cdot\Gamma)|^2})\big|\lesssim |\nabla_u(\xi\cdot\Gamma)|^{-1}\big(\lambda^{-\epsilon_1}|\nabla_u(\xi\cdot\Gamma)|\big)^{|\alpha|}
    \Ee
    on $\supp a_{\lambda,1}$.
    Considering that
    \[
    |(D^\intercal)^N a_{\lambda,1}|\lesssim \sum_{|\alpha^1|+\cdots|\alpha^{N+1}|= N}|\partial_u^{\alpha^{N+1}}a_{\lambda,1}|\prod_{n=1}^N|\partial_u^{\alpha^n}(\frac{\nabla_u(\xi\cdot\Gamma)}{|\nabla_u(\xi\cdot\Gamma)|^2})|,
    \]
    the inequalities \eqref{ineq:paa1} and \eqref{ineq:pgamma} imply that the desired estimate \eqref{upbd:ae1} holds.

\subsection*{When $p=\infty$}
    Here, we use kernel estimates for $\cA[a_\lambda]=\cA[a]-\cA[a\chi_{\lambda}]$. 
     Define the kernel $K[a]$ as follows:
    \[
    K[a](t,u,x):=(2\pi)^{-d-1}\iiint e^{i(t-t',x-t'\Gamma(u))\cdot(\tau,\xi)}a(t',u,\tau,\xi)d\xi d\tau dt'.
    \]
    Then  we can write $\mathcal A[a_\lambda]f(t,x)=\int K[a_\lambda](t,u,\cdot)\ast f(x)du$. 
    By Young's convolution inequality and the support condition $\supp_u a_\lambda\subset B(0,1)$, to show \eqref{aa01} for $p=\infty$,  it is enough to show that 
    \[
    \int |K[a_\lambda](t,u,x)|dx \lesssim 1,
    \]
    for all $(t,u)\in \R\times B(0,1)$.
    
    In the rest of this section, fix $(t,u)\in \R\times B(0,1)$.
    Since $a_\lambda=a-a\chi_\lambda$, the matter is  reduced to proving the following two inequalities: 
    \[
    \int |K[a](t,u,x)|dx \lesssim 1 \text{ and }\int |K[a\chi_\lambda](t,u,x)|dx \lesssim 1.
    \]
 The first one
    is easily verified by \eqref{cond:a}.
    Indeed, using change of variables $(\tau,\xi)\mapsto \lambda(\tau,\xi)$, and then integration by parts with respect to $(\tau,\xi)$-variables, we have
    \[
        |K[a](t,u,x)|\lesssim \lambda^{d+1} \int_{1/2}^4 (1+\lambda|(t-t',x-t'\Gamma(u))|)^{-d-2}dt'.
    \]
    Here, we also use $\supp_{\tau,\xi}a\subset \mathbb A_\lambda$.
 The above inequality implies the first inequality.

 Now,  we only need to establish the estimate $\|K[a\chi_\lambda](t,u,\cdot)\|_{L^1}\lesssim 1$.
    To achieve this, we use a similar approach as before, but with a different change of variables.
For this purpose, we  define a linear map $\mathcal L_u^\lambda:\R^{d}\rightarrow\R^{d}$ as follow: 
    \[
    \mathcal L_u^\lambda\partial_j\Gamma(u)=\lambda^{1/2+\epsilon_1}\partial_j\Gamma(u), \quad  \text{ for }j=1,\cdots,d-1,
    \]
    and
    \[
    \mathcal L_u^\lambda (0,\cdots,0,1)=\lambda (0,\cdots,0,1).
    \]
     And define a linear map $L_u^\lambda:\R^{d+1}\rightarrow\R^{d+1}$ by
 \[
 L_u^\lambda(\tau,\xi):=(\lambda^{\epsilon_1}\tau-\Gamma(u)\cdot (\mathcal L_u^\lambda)^\intercal \xi,(\mathcal L_u^\lambda)^\intercal\xi).
 \]
 Here, by $(\mathcal L_u^\lambda)^\intercal$ we denote the transpose of the linear operator $\mathcal L_u^\lambda$.
    By definition, $L_u^\lambda$ satisfies the following properties:
    \[
    L_u^\lambda(\tau,\xi)\cdot(1,\Gamma(u))=\lambda^{\epsilon_1}\tau,
    \]
    \[
     L_u^\lambda(\tau, \xi) \cdot(0,\partial_j\Gamma(u)) =\lambda^{1/2+\epsilon_1}\xi\cdot\partial_j\Gamma(u),
    \]
    for $j=1,\cdots,d-1$, and
    \[
      L_u^\lambda (\tau,\xi) \cdot (0,0,\cdots,0,1)=\lambda \xi\cdot (0,\cdots,0,1).
    \]
 From these properties and \eqref{maincaes}, we can get the following proposition, which says that after applying the change of variables $L_u^\lambda$, $a\chi_\lambda$ is ``well-normalized".
 \begin{prop}\label{prop:atu}
     Let $(t',u)\in [1/2,4]\times B(0,1)$. Set $ a_{t',u}^\lambda(\tau,\xi)= (a\chi_\lambda)(t',u, L_u^\lambda(\tau,\xi))$. Then
     \Be
\label{suppa0L}
|\supp  a_{t',u}^\lambda|\lesssim 1
\Ee
 and for $|\beta|\le d+2$,
\Be
\label{ineq:a0L}
|\partial_{(\tau,\xi)}^\beta  a_{t',u}^\lambda|\lesssim 1.
\Ee
 \end{prop}
 \begin{proof}
      We first prove the support condition \eqref{suppa0L}. Assume that $(\tau,\xi)\in \supp a_{t',u}^\lambda$. Since $(t',u,L_u^\lambda(\tau,\xi))\in \supp(a\chi_\lambda)$, $|L_u^\lambda(\tau,\xi)|\le2\lambda$ and \eqref{maincaes} tells that 
      \[
      |L_u^\lambda(\tau,\xi)\cdot\partial_u^\alpha(1,\Gamma(u))|\le 2\lambda^{\epsilon_1+|\alpha|/2}
      \]
      for $|\alpha|\le1$. By this and the properties of $L_u^\lambda$, we can see that 
      \Be
      \label{ineq:Mu}
      |M(u)\cdot (\tau,\xi)|\le 2d+2,
      \Ee
      where $M(u)$ is the $(d+1)\times(d+1)$ matrix given by
      \[
      M(u)=\begin{bmatrix}
1 &0& \cdots & 0& 0  \\
0 &1 & \ddots & \vdots& \partial_1\gamma(u) \\
\vdots & \ddots& \ddots &0 &\vdots \\
\vdots & &\ddots & 1& \partial_{d-1}\gamma(u) \\
0 & \cdots&\cdots & 0 & 1
\end{bmatrix}.
      \]
      Indeed, one can observe that $M(u)\cdot(\tau,\xi)=(\tau,\partial_1\Gamma(u)\cdot \xi,\cdots,\partial_{d-1}\Gamma(u)\cdot\xi,\xi_d)$, where the absolute value of each component is bounded above by $2$ due to   $(t',u,L_u^\lambda(\tau,\xi))\in \supp a\chi_\lambda$.
      
  By \eqref{ineq:Mu}, we can get $(\tau,\xi)\in M(u)^{-1}(B_{d+1}(0,2d+2))$ where $B_{n}(0,r):=\{x\in \R^n:|x|\le r\}$. Thus, we obtain that
     \[
     \supp a_{t',u}^\lambda\subset M(u)^{-1}(B_{d+1}(0,2d+2)).
     \]
     Since $M(u)$ is an upper triangular matrix with all diagonal entries equal to $1$, one can check that $\det(M(u))=1$, which implies  the desired \eqref{suppa0L}.

     Now, we show estimates for the derivatives \eqref{ineq:a0L}. By the properties of $L_u^\lambda$, we have that
     \[
a_{t',u}^\lambda(\tau,\xi)=a(t',u,L_u^\lambda(\tau,\xi))\eta_0(\tau)\prod_{i=1}^d\eta_0(\xi\cdot \partial_j\Gamma(u)).
\]
Thus, \eqref{ineq:a0L} follows from
 \eqref{cond:a} and the fact that $\|L_u^\lambda\|_{op}\lesssim \lambda$ where $\|L_u^\lambda\|_{op}=\sup_{(\tau,\xi)}|L_u^\lambda(\tau,\xi)|/|(\tau,\xi)|$.
 \end{proof}

Using Proposition \ref{prop:atu}, we can  prove the following lemma, which implies directly the desired estimate, $\|K[a\chi_\lambda](t,u,\cdot)\|_{L^1}\lesssim 1$.
 \begin{lem}\label{lem:ker}
        Let $g:\R\times \R^{d-1}\rightarrow \R$ be a function satisfying $|g|\le 1$. And set $\mathfrak a(t',u,\tau,\xi)=g(t',u)(a\chi_\lambda)(t',u,\tau,\xi)$. Then 
        \[
        \int |K[\mathfrak a](t,u,x)| dx\lesssim 1.
        \]
    \end{lem}
 
    \begin{proof}
        
    Using change of variables $(\tau,\xi)\mapsto   L_u^\lambda(\tau,\xi)$, we can get
    \[
    K[\mathfrak a](t,u,x)=C_d(\lambda)\iiint e^{i(t-t',x-t\Gamma(u))\cdot(\lambda^{\epsilon_1}\tau,\mathcal L_u^\lambda\xi)}\mathfrak a(t',u,  L_u^\lambda(\tau,\xi))d\xi d\tau dt',
    \]
where $C_d(\lambda)=(2\pi)^{-d-1}\lambda^{(d+1)/{2}+d\epsilon_1}$. By proposition \ref{prop:atu}, the symbol $\mathfrak a$ satisfies the followings: $|\supp \mathfrak a(t',u,L_u^\lambda(\cdot))|\lesssim 1$ and $|\partial_{(\tau,\xi)}^\beta(\mathfrak a(t',u,L_u^\lambda(\tau,\xi)))|\lesssim 1$ for $|\beta|\le d+2$.

Using integration by parts with respect to $(\tau,\xi)$-variables, we finally obtain
\[
\int |K[\mathfrak a](t,u,x)|dx\lesssim C_d(\lambda)\iint (1+|(\lambda^{\epsilon_1}(t-t'),(\mathcal L^\lambda_u)^\intercal(x-t\Gamma(u)))|^{-d-2}) dt'dx.
\]
Since $\lambda^{\epsilon_1}|\det \mathcal L^\lambda_u|\sim C_d(\lambda) $, we can conclude that the desired estimate holds.
    \end{proof}

\section{Decomposition away from degeneracy}\label{sec:decomp}
In this section, we prove \eqref{aa00}. To achieve this, we decompose the symbol based on degeneracy of $\gamma$. Specifically, we decompose the symbol with respect to the second derivatives, $(\sum_{|\alpha|=2}|\partial^\alpha\gamma(u)|^2)^{1/2}$. Additionally, our analysis requires a suitable upper bound on higher-order derivatives. For this reason, we further decompose symbols using
\[
\mathfrak D^m\gamma(u):=(\sum_{|\alpha|=m}|\partial^\alpha\gamma(u)|^2)^{1/2},
\]
for $m=2,\cdots,k$.

Before proceeding with the degeneracy decomposition, we first provide a basic decomposition.
Since we allow $\epsilon$-loss in \eqref{aa00}, it is enough to consider that $u$ is contained in a small ball.
Precisely, we decompose $a \chi_\lambda$ as follow:
\[
a \chi_\lambda=\sum_{\nu\in \Z^{d-1}}a\chi_\lambda\zeta(\lambda^{\epsilon_1}u-\nu)=:\sum_{\nu\in\Z^{d-1}} a_\nu,
\]
where $\zeta$ is a smooth function satisfying $\supp \zeta\subset B(0,1)$ and $\sum_{\nu\in\Z^{d-1}}\zeta(u-\nu)\equiv 1$. 
By triangle inequality, \eqref{aa00} is reduced to show that there exists  a constant  $C$ depending only on $d$ such that
\Be
    \label{ineq:main}
    \|\mathcal A[a_\nu]f\|_{L^p(\R^{d+1})}\lesssim  \lambda^{-\frac{2}p+C\epsilon_1}\|f\|_p,
    \Ee
    for all $\nu$. Indeed, \eqref{ineq:main} implies that \eqref{aa00} after choosing $C_\ast=C+d-1$.

    In the rest of this section, we show \eqref{ineq:main} for fixed $\nu$.
To do further decomposition with respect to degeneracy,
define $\eta_1=\eta_0(2^{-1}\cdot)-\eta_0(\cdot)$. Note that $\supp\eta_1\subset [-4,-1]\cup[1,4]$ and $\eta_0(2^{-1}\cdot)+\sum_{j\ge1}\eta_1(2^{-j}\cdot)\equiv1$.
We can observe that
\[
\prod_{m=2}^{k}\big(\eta_0(2^{-1}\lambda^{\frac 12}\mathfrak D^m\gamma(u))+\sum_{\ell_m\ge1} \eta_1(\lambda^{\frac 12}2^{-\ell_m}\mathfrak D^m\gamma(u))\big)\equiv 1.
\]
Set for $\vec\ell:=(\ell_2,\cdots,\ell_{k})$,
\[
\mathfrak a_{\vec\ell}:=a_\nu\prod_{m=2}^k\eta_1(\lambda^{1/2}2^{-\ell_m}\mathfrak D^m\gamma(u)) \text{ and } \mathfrak a_{0}:=a_\nu-\sum_{\vec\ell\in \N^{k-1}}\mathfrak a_{\vec\ell}.
\]
    By the assumption \eqref{cond:gamma} on $\gamma$, 
    we first note that $\mathfrak a_{\vec\ell}\equiv0$ unless the following conditions are satisfied:
\Be
\label{cond:lm}
2^{\ell_m}\lesssim \lambda^{\frac 12},
\Ee
for all $m=2,\cdots,k-1$ and
\Be
\label{cond:lk}
2^{\ell_k}\sim \lambda^{\frac 12}.
\Ee
Furthermore, when $k=2$, one can easily check  that $\mathfrak a_0\equiv 0$ since $a_\nu=\sum_{\ell_2\sim \lambda^{1/2}}\mathfrak a_\vl$.

By definition, we can get 
\Be
\label{eq:maindecomp}
\mathcal A[a_\nu]f=\mathcal A[\mathfrak a_{0}]f+\sum_{\vec\ell\in\mathbb L}\mathcal A[\mathfrak a_{\vec\ell}]f,
\Ee
where $\mathbb L:=\{\vl=(\ell_2,\cdots,\ell_k)\in \N^{k-1}:\eqref{cond:lm} \text{ and } \eqref{cond:lk} \text{ hold}  \}$.
We note that the number of elements in $\mathbb L$ satisfies that $ |\mathbb L|\lesssim \lambda^{\epsilon_1}$. Thus, to show \eqref{ineq:main}, it is sufficient to verify the following two propositions:
\begin{prop}\label{prop:aa0}
    Let $2\le p\le\infty$ and $k\ge3$. Then 
    \[
    \|\mathcal A[\mathfrak a_{0}]f\|_p\lesssim_{\epsilon_1} \lambda^{-\frac1{2(k-2)}+\frac{\epsilon_1}2}\|f\|_p.
    \]
\end{prop}
\begin{prop}\label{prop:main}
Let $p\ge\max\{6,4k-4\}$ and $\vl\in\mathbb L$. Then there exists a constant $C$ depending only on $d$ such that
\[
    \|\cA[\mathfrak a_\vl]f\|_p\lesssim \lambda^{C\epsilon_1-\frac 2p}\|f\|_p.
\]
\end{prop}
In our proofs of Proposition \ref{prop:aa0} and Proposition \ref{prop:main}, we crucially utilize the fact that the sizes of supports $|\supp_u\mathfrak a_0|$ and $|\supp_u\mathfrak a_\vl|$ are sufficiently small. Indeed, Proposition \ref{prop:aa0} is a direct consequence of the size estimate \eqref{phi-1} below. In the next subsection, we give a proof of Proposition \ref{prop:aa0} using an estimate on the support of $\mathfrak a_0$.

\subsection{Proof of Proposition \ref{prop:aa0}}\label{sec:a0}
In this subsection, we assume that $k\ge3$. 
We begin by proving  Proposition \ref{prop:aa0}, under the assumption that the following inequality holds:
\begin{equation}
\label{phi-1}
|\supp_u \mathfrak a_{0}|\lesssim_{\epsilon_1} \lambda^{-\frac{1}{2(k-2)}}.
\end{equation}
To show Proposition \ref{prop:aa0}, it suffices to show the $L^2$  and $L^\infty$ estimate since  interpolation between them provides $L^p$ estimates for $2\le p\le \infty$. 

When $p=2$, by  Plancherel's theorem, we have 
    \[
    \|\cA[\mathfrak a_0]f\|_2^2\le
    \big(\sup_{\xi}\int(\iint |\mathfrak a_{0}(t',u,\tau,\xi)|dudt')^2d\tau\big)\|f\|_2^2.
    \]
    For estimating the first term on the right-hand side,
 we apply
     Minkowski's inequality, yielding that
    \[
    \begin{aligned}
         \int(\iint |\mathfrak a_{0}(t',u,\tau,\xi)|dudt')^2d\tau  \lesssim  \big(\iint  (\int|\mathfrak a_{0}(t',u,\tau,\xi)|^2 d\tau)^{\frac 12} dudt'\big)^2. 
    \end{aligned}
    \]
    Since \eqref{maincaes} holds on $\supp \mathfrak a_0$, we can get $|\supp \mathfrak a_0(t',u,\cdot,\xi)|\lesssim \lambda^{\epsilon_1}$. Thus,  \eqref{phi-1}  gives that $\int(\iint |\mathfrak a_{0}(t',u,\tau,\xi)|dudt')^2d\tau\lesssim \lambda^{-1/(k-2)+\epsilon_1}$, which implies the desired $L^2$ estimate,
    \[
    \|\cA[\mathfrak a_0]f\|_2\lesssim \lambda^{-\frac{1}{2(k-2)}+\frac {\epsilon_1}2}\|f\|_2.
    \]
   
    When $p=\infty$, we use the fact that $\mathfrak a_0(t',u,\tau,\xi)=g(t',u)(a\chi_\lambda)(t',u,\tau,\xi)$ for a function $g$ satisfying $|g|\le1$.
    Then, Lemma \ref{lem:ker} implies that 
 $\|K[\mathfrak a_{0}]
    (t,u,\cdot)\|_{L^1}\lesssim 1$.
    Combined with $|\supp_u K[\mathfrak a_0]|\lesssim \lambda^{-1/(2(k-2))}$, which is a direct consequence of  \eqref{phi-1}, we can obtain the desired $L^\infty$ estimate,
    \[\|\mathcal A[\mathfrak a_{0}]f\|_\infty \lesssim   \lambda^{-1/2(k-2)}\|f\|_\infty.\]

Now, it remains to prove \eqref{phi-1}. To show this, we first check useful properties of $\mathfrak a_0$.
By the definition of $\mathfrak a_0$, there exists $m'=m'(u)\in\{2,\cdots, k-1\}$ such that
\Be
\label{dm'gamma}
\mathfrak D^{m'}\gamma(u)\le 4\lambda^{-1/2},
\Ee
for $u\in \supp_u \mathfrak a_0$. Using \eqref{dm'gamma} and $\mathfrak D^k\gamma(u)\sim 1$, we can observe that there exists $m=m(u)\in \{2,\cdots, k-1\}$ such that 
\begin{equation}
    \label{phi-1j}
    \mathfrak D^{m+1}\gamma(u)\ge \tilde {\mathcal B}^{-1} \lambda^{-\frac{k-m-1}{2(k-2)}},
    \quad \text{ and }\quad
    \mathfrak D^{m}\gamma(u)\le 8\lambda^{-\frac{k-m}{2(k-2)}},
\end{equation}
for $u\in \supp_u \mathfrak a_{0}$, 
where $\tilde {\mathcal B}\ge 1$ is a constant depending only on $ \mathcal B$ and $k$ satisfying $\mathfrak D^{k}\gamma(u)\ge \tilde {\mathcal B}^{-1}$ for all $u\in B(0,1)$.
Indeed, suppose that \eqref{phi-1j} fails for all $m\in \{2,\cdots,k-1\}$ for some $u\in \supp_u \mathfrak a_0$.
Then, since \eqref{phi-1j} fails for $m=k-1$, we have $\mathfrak D^{k-1}\gamma(u)\ge 8 \lambda^{-1/2(k-2)}$. Inductively, we can see  
 that $\mathfrak D^{m}\gamma(u)\ge 8 \lambda^{(m-k)/2(k-2)}\ge 8\lambda^{-1/2}$ for all $m\in \{2,\cdots, k-1\}$, which gives a contradiction considering \eqref{dm'gamma}.

Thus,  $\supp_u \mathfrak a_0$  is contained in $\cup_{m=2}^{k-1}S_{m}$ , where
\[
S_{m}=\{u\in B(\lambda^{-\epsilon_1}\nu,\lambda^{-\epsilon_1}): \eqref{phi-1j} \text{ holds}\}.
\]
For showing \eqref{phi-1},  it is sufficient to show that 
\[
|S_m|\lesssim \lambda^{-\frac 1{2(k-2)}}.
\]
It is a direct consequence of the
 following lemma:
\begin{lem}\label{lem:S}
    Let $m=2,\cdots,k-1$, and $A,B\ge\lambda^{-1/2}$. Assume that
    $S\subset B(\lambda^{-\epsilon_1}\nu,\lambda^{-\epsilon_1})$ satisfies $\mathfrak D^{m+1}\gamma(u)\ge A$ and $\mathfrak D^{m}\gamma(u)\le B$ for  $u\in S$. Then $|S|\lesssim  A^{-1}B$.
\end{lem}
\begin{proof}
    We can observe that 
$
|S|\le \int_{S} A^{-1}\mathfrak D^{m+1}\gamma(u) du.
$
Since $\mathfrak D^{m}\gamma\lesssim \sum_{|\alpha|=m}|\partial^\alpha\gamma|$, it is sufficient to show that for $\alpha\in \N_0^{d-1}$ satisfying $|\alpha|=m+1$,
\Be
\nonumber
\int_S|\partial^\alpha\gamma(u)|du \lesssim B.
\Ee
For proving this, we use the Taylor expansion of $\gamma$.
Set $p_\nu$ is the $N_{\epsilon_1}$-th order Taylor polynomial of $\gamma$ at the point $\lambda^{-\epsilon_1}\nu$ for sufficiently large $N_{\epsilon_1}>k+1/\epsilon_1$.
By Taylor's theorem, it holds that
\Be
\label{pnu}
\sup_{|\alpha|\le k}|\partial^\alpha \gamma(u)- \partial^\alpha p_{\nu}(u)|\lesssim\lambda^{-1},
\Ee 
for $u\in B(\lambda^{-\epsilon_1}\nu,\lambda^{-\epsilon_1})$.  By triangle inequality with \eqref{pnu} and the condition $\lambda^{-1}\le B$, it is enough to show that for $|\alpha|=m+1$,
\Be
\label{ineq:gammaB}
\int_S|\partial^\alpha p_\nu(u)|du\lesssim B.
\Ee

To show \eqref{ineq:gammaB}, fix $\alpha=(\alpha_1,\cdots,\alpha_{d-1})$. And choose $j$ satisfying $\alpha_j\ge1$.
Since $p_\nu$ is a polynomial,  
 there exist intervals $I_n^\alpha=I_n^\alpha(u_1,\cdots,u_{j-1},u_{j+1},\cdots,u_{d-1})$ and a integer $N'\le N_{\epsilon_1}+1$ such that $\cup_{n=1}^{N'}I_n^\alpha=\R$ and  $u_j\mapsto\partial^\alpha p_\nu(u)$ is nonnegative or nonpositive on each $I_n^\alpha$. 
By this observation, we can check that
\[
\int_{\{u_j:u\in S\}}|\partial^\alpha p_\nu(u)|du_j\le \sum_{n=1}^{N'}\big|\int_{\{u_j\in I_n:u\in S\}} \partial^\alpha p_\nu(u) du_j\big|.
\]
By the Fundamental theorem of calculus, we have that
the right-hand side of the above inequality is bounded by 
$2N' \max_{u\in S} |\partial^{\alpha'} p_\nu (u)|$ where $\alpha'=(\alpha_1,\cdots,\alpha_{j-1},\alpha_j-1,\alpha_{j+1},\cdots,\alpha_{d-1})$. Since $|\alpha'|=|\alpha|-1=m$, \eqref{pnu} and the condition $\mathfrak D^m\gamma(u)\le B$ tells that
\[
\int_{\{u_j:u\in S\}}|\partial^\alpha p_\nu(u)|du_j\lesssim B. 
\]
Thus, Fubini's theorem implies the desired estimate \eqref{ineq:gammaB}. 
\end{proof}

 In the rest of this section, we give a proof of Proposition \ref{prop:main}. To do this, we prove some properties of $\mathfrak a_\vl$ in the next subsection.


\subsection{Properties of $\mathfrak a_\vl$}

Let $\vl\in \mathbb L$.
     Define a set 
     \[
     S_{\vec{\ell}}:=\{u\in B(\lambda^{-\epsilon_1}\nu,\lambda^{-\epsilon_1}): 1\le \lambda^{1/2} 2^{-\ell_m}\mathfrak D^m\gamma(u)\le 4 
     \text{ for } m=2,\cdots,k\}.
     \]
     Note that $\supp_u \mathfrak a_{\vec\ell}\subset S_{\vec\ell}$.
     By Lemma \ref{lem:S}, the set $S_{\vec\ell}$ satisfies that $|S_{\vec\ell}|\lesssim \delta_{\vec\ell}$ where
     \[
     \delta_{\vec{\ell}}:=\begin{cases}
         \min_{2\le m\le k-1} 2^{\ell_m-\ell_{m+1}}, & \text{if } k\ge 3,\\
         1, &\text{if } k=2.
     \end{cases}
     \]
      To prove Proposition \ref{prop:main}, we require an additional property of the support of $\mathfrak a_\vl$.  In this subsection, we demonstrate that  $S_\vl$ is covered by $O(\delta_\vl^{-d+2})$-many balls with radius $\delta_\vl$. This property simplifies the analysis of $\mathfrak a_\vl$ in the subsequent sections by allowing us to focus on  these $\delta_\vl$-balls.

     To establish the desired property, we first show that $\mathfrak D^m\gamma$ has a suitable bound on $\delta_\vl$-neighborhood of  $S_\vl$.
     Since $|\nabla\mathfrak D^m\gamma|\lesssim \mathfrak D^{m+1}\gamma$, we can verify the following proposition.
     \begin{prop}\label{prop:Dmg}
         Let $\vec\ell\in \mathbb L$ and $m=2,\cdots,k$. Suppose $u_0\in S_{\vec\ell}$. Then we have
          for $u\in B(u_0,2\delta_{\vec\ell})$, 
             \Be
          \nonumber
             \mathfrak D^{m}\gamma(u)\lesssim \lambda^{-1/2}2^{\ell_m}.
             \Ee
     \end{prop}
     \begin{proof}
      We prove this by induction. By \eqref{cond:lk}, the above inequality holds for $m=k$. Now, for $2\le m\le k-1$, suppose that $             \mathfrak D^{m+1}\gamma(u)\lesssim \lambda^{-1/2}2^{\ell_{m+1}}$. Then Mean value theorem implies that
      \[
      \mathfrak D^{m}\gamma(u)\lesssim \mathfrak D^{m}\gamma(u_0)+\lambda^{-1/2}2^{\ell_{m+1}}\delta_{\vec\ell}.
      \]
      Here, we use the fact that $|\nabla\mathfrak D^{m}\gamma|\lesssim \mathfrak D^{m+1}\gamma$.
It gives that $             \mathfrak D^{m}\gamma(u)\lesssim \lambda^{-1/2}2^{\ell_m}$, since $\mathfrak D^{m}\gamma(u_0)\lesssim \lambda^{-1/2}2^{\ell_{m}}$ and $2^{\ell_{m+1}}\delta_{\vec\ell}\le 2^{\ell_{m}}$. Thus, we can conclude the desired inequality for all $m=2,\cdots, k$.
     \end{proof}
For later use, we choose a constant $\mathfrak C\sim 1$ such that
             \Be
             \label{cond:Dmg2}
             \mathfrak D^{m}\gamma(u)\le\frac{\mathfrak C}{100d(k+1)^d} \lambda^{-1/2}2^{\ell_m},
             \Ee
             for all $m=2,\cdots,k$ and $u$ in the $\delta_\vl$-neighborhood of $S_\vl$.
     This condition allows us to effectively control the derivatives of $\gamma$ on the $(\mathfrak C^{-1}\delta_\vl)$-neighborhood of $S_\vl$. Using this, we can  estimate the number of elements in the set
     \[
     \mathcal I_{\vec \ell}:=\{z\in \delta_{\vec{\ell}}\Z^{d-1}: B(z,\delta_{\vec{\ell}})\cap S_{\vec{\ell}}\neq \emptyset\}.
     \]
     Since $|S_\vl|\lesssim \delta_\vl$ and $|B(z,\delta_\vl)|\lesssim \delta_\vl^{d-1}$, one can see that $ |\mathcal I_\vl|\gtrsim \delta_\vl^{-d+2}$.
     Further, 
     $\mathcal I_{\vec{\ell}}$ satisfies the following proposition, which says that $S_\vl$ can be covered by $O(\delta_\vl^{-d+2})$-many balls $\{B(z,\delta_{\vec{\ell}})\}_{z\in \mathcal I_\vl}$.
\begin{prop}\label{prop:il}
Let $\vl\in\mathbb L$. Then
$ |\mathcal I_{\vec\ell}|\lesssim \delta_{\vec{\ell}}^{-(d-2)}$.
\end{prop}

\begin{proof}
First, note that
  by \eqref{cond:lm} and \eqref{cond:lk} one can check that
     \Be
\nonumber
     \delta_\vl\lesssim 1.
     \Ee
     
Define a set 
\[
\tilde S_{\vec \ell}:=\{u\in B(\lambda^{-\epsilon_1}\nu,\lambda^{-\epsilon}): 1/2\le \lambda^{1/2} 2^{-\ell_m}\mathfrak D^m\gamma(u)\le 8 
     \text{ for } m=2,\cdots,k\}.
     \]
     Since $\{B(z,\delta_\vl)\}_{z\in \mathcal I_\vl}$ are finitely overlapping, we can get
     \[
     |\tilde S_\vl|\gtrsim \sum_{z\in \mathcal I_\vl} |B(z,\delta_\vl)\cap \tilde S_\vl|.
     \]
     From the observation $|\tilde S_\vl|\lesssim\delta_\vl$ by Lemma \ref{lem:S}, proving Proposition \ref{prop:il} is reduced to showing that for all $z\in \mathcal I_\vl$,
     \Be
     \nonumber
     |B(z,\delta_\vl)\cap \tilde S_\vl|\gtrsim \delta_\vl^{d-1}.
     \Ee

We use \eqref{cond:Dmg2} to obtain this.
     For $z\in \mathcal I_\vl$, we can choose a point $u(z)\in S_\vl\cap B(z,\delta_\vl)$. 
     By Mean value theorem and the fact that $|\nabla\mathfrak D^{m}\gamma|\le (m+1)^d \mathfrak D^{m+1}\gamma$, 
        \eqref{cond:Dmg2} tells that 
     $B(u(z),\mathfrak C^{-1}\delta_\vl)\subset \tilde S_\vl$. Thus, it is enough to verify that $|B(z,\delta_\vl)\cap B(u(z),\mathfrak C^{-1}\delta_\vl)|\gtrsim \delta_\vl^{d-1}$, which is a direct consequence of the facts that $u(z)$ is contained in $ B(z,\delta_\vl)$ and $\mathfrak C\sim 1$.
\end{proof}

\subsection{Proof of Proposition \ref{prop:main}}

In this subsection, we verify Proposition \ref{prop:main}. 
We first treat the easier case, when $ \delta_\vl \le \mathfrak C\lambda^{\epsilon_1}2^{-\ell_2}$. Note that this case can not arise if $k=2$, since $2^{\ell_2}\sim \lambda^{1/2}$. Recall that $|\supp_u \mathfrak a_\vl|\lesssim \delta_\vl$. By the same argument in Section \ref{sec:a0}, we can deduce that \[\|\cA [\mathfrak a_\vl]f\|_p\lesssim \delta_\vl \lambda^{\frac{\epsilon_1}{2}}\|f\|_p.\] Proposition \ref{prop:main} for this case follows if we show that $\delta_\vl\lesssim \lambda^{-1/(2k-2)+\epsilon_1}$. For this, since $\delta_\vl^{k-2}\le 2^{\ell_2-\ell_k}\sim \lambda^{-1/2}2^{\ell_2}$, we can observe that
\Be
\label{ineq:deltal}
\delta_\vl\lesssim (\lambda^{-\frac 12}2^{\ell_2})^{\frac 1{k-2}},
\Ee
when $k\ge 3$.
Considering the fact that $\min\{\lambda^{\epsilon_1}2^{-\ell_2},(\lambda^{-\frac 12}2^{\ell_2})^{\frac 1{k-2}}\}\le \lambda^{-1/(2k-2)+\epsilon_1}$,  we can conclude the desired,  $\delta_\vl\lesssim \lambda^{-1/(2k-2)+\epsilon_1}$.

Now, we consider the main case, when $\lambda^{\epsilon_1}2^{-\ell_2}\le \mathfrak C^{-1}\delta_\vl$. For this case, we use decoupling inequality. Set 
\[
\mathbb L_0=\{\vl\in \mathbb L:\lambda^{\epsilon_1}2^{-\ell_2}\le \mathfrak C^{-1}\delta_\vl\}.
\]
For $\vl\in \mathbb L_0$, $\delta>0$, and $\mu\in \delta \Z^{d-1}$, we define
\[
\mathfrak a_{\vec\ell,\mu}^\delta:=\mathfrak a_{\vec\ell}\,\zeta(\delta^{-1}(u-\mu)).
\]
We can decompose $a_\vl$ as follow:
$
\mathfrak a_\vl=\sum_{z\in \mathcal I_\vl} \mathfrak a_{\vl,z}^{\delta_\vl}.
$
By H\"older's inequality combined with Proposition \ref{prop:il}, we can get 
\Be
\label{trivialdecoupling}
\|\cA[\mathfrak a_\vl]f\|_p\lesssim \delta_\vl^{-(d-2)(1-\frac 1p)}\big(\sum_{z\in \mathcal I_\vl}\|\cA[\mathfrak a_{\vl,z}^{\delta_\vl}]f\|_p^p\big)^{\frac 1p}.
\Ee
The proof of Proposition \ref{prop:main} is based on the following decoupling inequality. 

\begin{thm}\label{thm:avl}
    Let $\vl\in \mathbb L_0$ and $p\ge 6$. Then there exists a constant $C>0$ depending only on $d$ such that
\[
\begin{aligned}
    (\sum_{z\in \mathcal I_\vl}\|\mathcal A[\mathfrak a_{\vec\ell,z}^{\delta_\vl}]f\|_p^p)^{\frac 1p}&\lesssim \lambda^{C\epsilon_1}(\delta_{\vec\ell}2^{\ell_2})^{d-1-\frac{d+2}{p}}(\sum_{\tilde z\in \tilde{\delta}_\vl \Z^{d-1}}\|\mathcal A[\mathfrak a_{\vec\ell,\tilde z}^{\tilde{\delta}_\vl}]f\|_p^p)^{\frac 1p},
\end{aligned}
\]
where $\tilde{\delta}_\vl:=\lambda^{\epsilon_1}2^{-\ell_2}$.
\end{thm}
By the above inequality, the symbol $\mathfrak a_{\vl,z}^{\delta_\vl}$ can be decomposed into smaller components $\mathfrak a_{\vec\ell,\tilde z}^{\tilde{\delta}_\vl}$. Then we can use smallness of the support $\mathfrak a_{\vec\ell,\tilde z}^{\tilde{\delta}_\vl}$. Utilizing support conditions, we can have the following lemma.
\begin{lem}\label{lem:support}
    Let $\vl\in \mathbb L_0$. Then
    \Be
    \label{disjointness}
(\sum_{\tilde z\in \tilde{\delta}_\vl \Z^{d-1}}\|\mathcal A[\mathfrak a_{\vec\ell,\tilde z}^{\tilde{\delta}_\vl}]f\|_p^p)^{\frac 1p}\lesssim \lambda^{\epsilon_1}\tilde{\delta}_\vl^{d-1-\frac{d-2}{p}}\|f\|_p.
\Ee

\end{lem}

Proposition \ref{prop:main} follows from Theorem \ref{thm:avl} and Lemma \ref{lem:support}.
Indeed, by \eqref{trivialdecoupling}, Theorem \ref{thm:avl}, and Lemma \ref{lem:support}, we can check that there exists a constant $C>0$ depending only on $d$ such that for $p\ge 6$,
\[
\|\mathcal A[\mathfrak a_\vl]f\|_p\lesssim \lambda^{C\epsilon_1}\delta_\vl^{1-\frac 4p}2^{-\frac{4\ell_2}{p}}\|f\|_p.
\]
This directly establishes Proposition \ref{prop:main} for the case $k=2$.
To verify Proposition \ref{prop:main} for $k\ge3$, it suffices to show that 
\[\delta_\vl^{1-\frac 4p}2^{-\frac{4\ell_2}{p}}\lesssim \lambda^{-\frac 2p},
\]
for $p\ge 4k-4$. This follows from the observation that
\[
\delta_\vl^{1-\frac 4p}2^{-\frac{4\ell_2}{p}}\lesssim (\lambda^{-\frac 12}2^{\ell_2})^{\frac 1{k-2}(1-\frac 4p)}2^{-\frac{4\ell_2}{p}},
\]
from  \eqref{ineq:deltal}. Notably, the right-hand side of this inequality increases with $\ell_2$, given that $p\ge 4k-4$. Thus, the condition \eqref{cond:lm} ensures the desired result.

Now, it remains to verify Theorem \ref{thm:avl} and Lemma \ref{lem:support}.
In the rest of this subsection, we prove Lemma \ref{lem:support}.
A proof of Theorem \ref{thm:avl} is provided in Section \ref{sec:decoup}.

\begin{proof}[Proof of Lemma \ref{lem:support}]
We prove Lemma \ref{lem:support} by interpolation and it suffices to show \eqref{disjointness} for the cases $p=2$ and $p=\infty$.
    When $p=\infty$, we use the kernel estimate, Lemma \ref{lem:ker}. One can observe that $\mathfrak a_{\vl,\tilde z}^{\tilde{\delta}_\vl}=g_{\vl,\tilde z}a\chi_\lambda$ where
    \[
    g_{\vl,\tilde z}(t,u)=\zeta(\lambda^{\epsilon_1}u-\nu)\zeta(\tilde{\delta}_\vl(u-\tilde z))\prod_{m=2}^k\eta_1(\lambda^{1/2}2^{-\ell_m}\mathfrak D^m\gamma(u)).
    \]
    By Lemma \ref{lem:ker}, we have that
    \[
    |\mathcal A[\mathfrak a_{\vec\ell,\tilde z}^{\tilde{\delta}_\vl}]f|\le \int_{B(\tilde z,\tilde{\delta}_\vl)}|K[\mathfrak a_{\vl,\tilde z}](t,u,\cdot)\ast f(x)|du\lesssim \tilde{\delta}_\vl^{d-1}\|f\|_\infty,
    \]
    which gives the desired \eqref{disjointness} for $p=\infty$.

For showing \eqref{disjointness} when $p=2$, we first investigate Fourier support of $\mathcal A[\mathfrak a_{\vec\ell,\tilde z}^{\tilde{\delta}_\vl}]f$. Since $\mathfrak a_{\vec\ell,\tilde z}^{\tilde{\delta}_\vl}\equiv 0$ if $S_\vl\cap B(\tilde z,\tilde{\delta}_\vl)=\emptyset$, we only need to consider $\tilde z$
contained in the set
\[
\tilde{\mathcal I}_\vl:=\{\tilde z\in \tilde{\delta}_\vl\Z^{d-1}:S_\vl\cap B(\tilde z,\tilde{\delta}_\vl)\neq\emptyset\}.
\]
By \eqref{maincaes}, if $\xi \in \supp_\xi \mathfrak a_{\vec\ell,\tilde z}^{\tilde{\delta}_\vl}$, then there exists a point $u\in B(\tilde z,\tilde \delta_\vl)$ such that $|\xi\cdot\partial_j\Gamma(u)|\le 2\lambda^{1/2+\epsilon_1}$ for $j=1,\cdots,d-1$.
Using Mean value theorem combined with the definition of $\tilde \delta_\vl$ and  Proposition \ref{prop:Dmg} for $m=2$,  we can observe that there exists a constant $C_0>1$ such that $|\xi\cdot\partial_j\Gamma(\tilde z)|\le C_0 \lambda^{1/2+\epsilon_1}$, for $\xi \in \supp_\xi \mathfrak a_{\vec\ell,\tilde z}^{\tilde{\delta}_\vl}$. 
Therefore, $\mathcal A[\mathfrak a_{\vec\ell,\tilde z}^{\tilde{\delta}_\vl}]f=\mathcal A[\mathfrak a_{\vec\ell,\tilde z}^{\tilde{\delta}_\vl}]P_{\tilde z}f$, where $P_{\tilde z}$ is frequency localization satisfying 
\[
\widehat{P_{\tilde z}f}(\xi)=\hat f(\xi)\prod_{j=1}^{d-1} \eta_0(C_0^{-1}\lambda^{-1/2-\epsilon_1}\xi\cdot \partial_j\Gamma(\tilde z)).
\]
Indeed, the Fourier support of $\mathcal A[\mathfrak a_{\vec\ell,\tilde z}^{\tilde{\delta}_\vl}]f$ is contained in $\R\times \supp \widehat{P_{\tilde z}f}$.

The Fourier transform of $P_{\tilde z}f$ is supported in a set $\mathcal S(\tilde z)$ where
\[
\mathcal S(\tilde z):=\{\xi\in  \mathbb A_\lambda: |\xi\cdot \partial_j\Gamma(\tilde z)|\le 2C_0 \lambda^{1/2+\epsilon_1} \text{ for }j=1,\cdots,d-1\}.
\]
Note that for $\xi\in \mathcal S(\tilde z)$, there exist $r\in \R$ and $v\in \R^d$ such that $\xi=r(-\nabla\gamma(\tilde z),1)+v$ where  $|r|\sim \lambda, |v|\lesssim  \lambda^{\frac 12+\epsilon_1}$. By this observation, we can see that there exists a constant $C_1>1$ such that if $|\nabla\gamma(\tilde z)-\nabla\gamma(\tilde z')|\ge C_1\lambda^{-\frac 12+\epsilon_1}$ then $\mathcal S(\tilde z)\cap \mathcal S(\tilde z')=\emptyset$. 
For a given $\tilde z\in \tilde{\mathcal I}_\vl$, define a set 
\[
\tilde{\mathcal I}_\vl(\tilde z):=\{\tilde z'\in \tilde{\mathcal I}_\vl:  |\nabla\gamma(\tilde z)-\nabla\gamma(\tilde z')|< C_1\lambda^{-\frac 12+\epsilon_1}\}.
\]
By definition, we can observe that if $\tilde z'\notin \tilde{\mathcal I}_\vl(\tilde z)$, then the Fourier supports of $P_{\tilde z}f$ and $P_{\tilde z'}f$ are disjoint. Thus, the number of the elements in the set $\tilde{\mathcal I}_\vl(\tilde z)$ gives an upper bound for the number of overlaps among the Fourier supports of $\{P_{\tilde z}f\}_{\tilde z}$. We claim that for $\tilde z\in \tilde{\mathcal I}_\vl$,
\Be
\label{ineq:overlaps}
 |\tilde{\mathcal I}_\vl(\tilde z)|\lesssim \tilde{\delta}_\vl^{-(d-2)}.
\Ee
Assuming \eqref{ineq:overlaps}, we prove \eqref{disjointness} for $p=2$.

By Plancherel's theorem, we have
\[
\begin{aligned}
\|\mathcal A[\mathfrak a_{\vec\ell,\tilde z}^{\tilde{\delta}_\vl}]f\|_2^2&\sim \iint (\iint e^{-it'(\tau+\Gamma(u)\cdot \xi)} \mathfrak a_{\vec\ell,\tilde z}^{\tilde{\delta}_\vl}(t',u,\tau,\xi)dudt'\big)^2|\widehat{P_{\tilde z} f}(\xi)|^2d\xi d\tau,    \\
&\lesssim \lambda^{\epsilon_1}\tilde{\delta}_\vl^{2d-2}\|P_{\tilde z}f\|_2^2.
\end{aligned}
\]
Here, we use Minkowski's inequality combined with the facts $|\supp \mathfrak a_{\vec\ell,\tilde z}^{\tilde{\delta}_\vl}(t',u,\cdot,\xi)|\lesssim \lambda^{\epsilon_1}$ and $\supp_u \mathfrak a_{\vec\ell,\tilde z}^{\tilde{\delta}_\vl}\subset B(\tilde z,\tilde{\delta}_\vl)$.
Using the above inequality, we can get 
\[
(\sum_{\tilde z\in \tilde{\mathcal I}_\vl}\|\mathcal A[\mathfrak a_{\vec\ell,\tilde z}^{\tilde{\delta}_\vl}]f\|_2^2)^{\frac 12}\lesssim \lambda^{\epsilon_1}\tilde{\delta}_\vl^{d-1}(\sum_{\tilde z\in \tilde {\mathcal I}_\vl}\|P_{\tilde z}f\|_2^2)^{\frac 12}.
\]
By \eqref{ineq:overlaps}, the Fourier supports of $P_{\tilde z}f$ overlap at most $O(\tilde{\delta}_\vl^{-(d-2)})$ times, which tells that 
\[
(\sum_{\tilde z\in \tilde {\mathcal I}_\vl}\|P_{\tilde z}f\|_2^2)^{\frac 12}\lesssim \tilde{\delta}_\vl^{-\frac{d-2}2}\|f\|_2.
\]
 As a consequence, we get the desired, \eqref{disjointness} for $p=2$.

Now, it remains to prove \eqref{ineq:overlaps}. By Proposition \ref{prop:il} and the observation that $\tilde{\mathcal I}_\vl\subset \cup_{z\in \mathcal I_\vl}B(z,\delta_\vl)$, \eqref{ineq:overlaps} follows from
\[
|\tilde{\mathcal I}_\vl(\tilde z)\cap B( z,\delta_\vl)|\lesssim (\delta_\vl/\tilde{\delta}_\vl)^{d-2},
\]
for $z\in \mathcal I_\vl$. We can decompose $B(z,\delta_\vl)$ into finite pieces $B(u,\mathfrak C^{-1}\delta_\vl)$. Thus, it is sufficient to verify that for any $u\in B(0,1)$
\Be
\label{ineq:localoverlaps}
|\tilde{\mathcal I}_\vl(\tilde z)\cap B( u,\mathfrak C^{-1}\delta_\vl)|\lesssim (\delta_\vl/\tilde{\delta}_\vl)^{d-2}.
\Ee

If $\tilde{\mathcal I}_\vl(\tilde z)\cap B( u,\mathfrak C^{-1}\delta_\vl)$ is not empty, we can easily check that $B(u,\mathfrak C^{-1}\delta_\vl)$ is contained in the $(3\mathfrak C^{-1}\delta_\vl)$-neighborhood of $S_\vl$. This follows from the facts that $\tilde{\mathcal I}_\vl(\tilde z)\subset \tilde{\mathcal I}_\vl$ and  $\tilde \delta_\vl\le \mathfrak C^{-1}\delta_\vl$.  By \eqref{cond:Dmg2}, the second order derivatives of $\gamma$ satisfy   $1/2\le\lambda^{1/2}2^{-\ell_2}\mathfrak D^2\gamma\le 8$ on $B(u,\mathfrak C^{-1}\delta_\vl)$. Furthermore, there exists a unit vector $v_{u}\in \R^{d-1}$ such that  
\[|\langle v_{u},\nabla\rangle^2\gamma| \ge \frac{1}{4d} \lambda^{-1/2}2^{\ell_2}\]   
on $B(u,\mathfrak C^{-1}\delta_\vl)$. 
This observation implies  the following:
\[ |\nabla\gamma(u_1 )-\nabla\gamma(u_2)|\ge2 C_1\lambda^{-1/2+\epsilon_1},\]
for any $u_1,u_2\in B(u,\mathfrak C^{-1}\delta_\vl)$ satisfying 
 $|\langle v_{u},u_1-u_2\rangle|\ge (8dC_1+20d)\lambda^{\epsilon_1}2^{-\ell_2}$ and $|\Pi_{v_{u}}(u_1-u_2)|\le\lambda^{\epsilon_1}2^{-\ell_2}$. Here, by $\Pi_{v}$ we denote the projection operator onto $\{v\}^\perp$. Thus, we conclude that if $\tilde z_1,\tilde z_2\in \tilde{\mathcal I}_\vl(\tilde z)\cap B( u,\mathfrak C^{-1}\delta_\vl)$ satisfy that $|\Pi_{v_{u}}(\tilde z_1-\tilde z_2)|\le\tilde \delta_\vl$, then $|\tilde z_1-\tilde z_2|\lesssim \tilde \delta_\vl$. This implies that for $u'\in B( u,\mathfrak C^{-1}\delta_\vl)$,
 \[
 |\{\tilde z'\in \tilde{\mathcal I}_\vl(\tilde z)\cap B( z',\mathfrak C^{-1}\delta_\vl):|\Pi_{v_{z'}}(\tilde z'-u')|\le\tilde \delta_\vl/2\}|\lesssim 1.
 \]
From this, we obtain the desired result, \eqref{ineq:localoverlaps}, since the ball $B(u,\mathfrak C^{-1}\delta_\vl)$ can be covered by $O((\delta_\vl/\tilde{\delta}_\vl)^{d-2})$-many tubes of the form $\{u'':|\Pi_{v_{u}}(u''-u')|\le\tilde \delta_\vl/2\}$.
\end{proof}

\section{Proof of Theorem \ref{thm:avl}}\label{sec:decoup}

In this section, we prove Theorem \ref{thm:avl} using decoupling inequality for cone in $\R^3$.
Let $\mathbf r:\R\rightarrow \R$ be a smooth function satisfying $\mathbf r''\neq 0$. Define a curve in $\R^{d+1}$, \[\tilde{\mathbf r}(s)=v_1+sv_2+\mathbf r(s)v_3,\]
where $\{v_i\}_{i=1}^3$ are orthonormal vectors in $\R^{d+1}$.
And by $\mathbf s(s,\delta;\tilde{\mathbf r})$ denote the set of $\omega\in \R^{d+1}$ which satisfies
\[
\begin{aligned}
    &1/2\le |\langle \tilde{\mathbf r}^{(2)}(s),\omega\rangle|\le 2,\\
    &|\langle\tilde{\mathbf r}^{(m)}(s),\omega\rangle|\le \delta^{2-m}, \quad m=0,1.
\end{aligned}
\]
Consider $\delta$-separated points $\{s_m\}_{m=1}^M$ satisfying $|s_m-s_{m'}|\ge\delta$ if $m\neq m'$ and $[0,1]\subset\cup_{m=1}^M(s_m-\delta,s_m+\delta)$. Then, we set
\[
\mathbf s_m=\mathbf s(s_m,\delta;\tilde{\mathbf r}).
\]

By Bourgain-Demeter \cite{BD} (also see \cite{PS}), the following theorem holds.
\begin{thm}[\cite{BD}]\label{BD decoup}
    Let $0<\delta<1$ and $ p\ge 6$. Suppose $\mathbf r(s)=cs^2$ for a constant $c\neq 0$. Then for $\epsilon>0$, there exists a constant $C_\epsilon$ depending only on $\epsilon$ and $c$ such that
    \[
    \|\sum_{1\le m\le M}f_m\|_p\le C_\epsilon \delta^{-1+\frac 4p-\epsilon}(\sum_{1\le m\le M}\|f_m\|_p^p)^{\frac 1p},
    \]
    whenever $\supp \hat{f}_m\subset \mathbf s_m$ for $1\le m\le M$.
\end{thm}
To establish Theorem \ref{thm:avl},  we use Proposition \ref{prop:avl} below, which allows us to decompose $\mathcal A[\mathfrak a_{\vl,\mu}^{\delta}]$ into slightly smaller components.
A key tool in  proving Proposition \ref{prop:avl} is  the decoupling inequality given in Theorem \ref{BD decoup}.

Let $ N_0$ denote the smallest integer satisfying  $\lambda^{-\epsilon_1\cdot ( 3/2)^{N_0}}\delta_{\vec\ell}\le \lambda^{\epsilon_1}2^{-\ell_2}$. By \eqref{cond:lm}, it follows that $N_0\lesssim \log(1/\epsilon_1)$.
Set 
\[
\delta_{\vec\ell,n}:=\lambda^{-\epsilon_1\cdot (\frac 32)^n}\delta_{\vec\ell} ,\quad \text{ for } 0\le n\le N_0-1,
\]
and $\delta_{\vl,N_0}:=\lambda^{\epsilon_1}2^{-\ell_2}$. 
For simplicity, we write $\delta_{n}:=\delta_{\vl,n}$. Note that
\[
\lambda^{\epsilon_1}2^{-\ell_2}=\delta_{N_0}\le\delta_{N_0-1}\le\cdots\le \delta_0=\lambda^{-\epsilon_1}\delta_\vl.
\]
Then, the following holds.
\begin{prop}\label{prop:avl}
    Let $0\le n\le N_0-1$ and $p\ge 6$. Then for $\epsilon>0$, we have
    \[
    (\sum_{\mu\in \delta_{ n}\Z^{d-1}}\|\mathcal A[\mathfrak a_{\vec\ell,\mu}^{\delta_{n}}]f\|_p^p)^{\frac 1p}\lesssim (\frac{\delta_{ n}}{\delta_{n+1}})^{d-1-\frac{d+2}{p}+\epsilon} (\sum_{\tilde \mu\in \delta_{n+1}\Z^{d-1}}\|\mathcal A[\mathfrak a_{\vec\ell,\tilde \mu}^{\delta_{n+1}}]f\|_p^p)^{\frac 1p}.
    \]
\end{prop}

Theorem \ref{thm:avl} follows from the repeated application of Proposition \ref{prop:avl}.

\begin{proof}[Proof of Theorem \ref{thm:avl} assuming Proposition \ref{prop:avl}]
After applying H\"older's inequality, we have
\[
(\sum_{z\in \mathcal I_\vl}\|\mathcal A[\mathfrak a_{\vec\ell,z}^{\delta_{\vec\ell}}]f\|_p^p)^{\frac 1p}\lesssim \lambda^{d\epsilon_1}(\sum_{\mu\in \delta_{0}\Z^{d-1}}\|\mathcal A[\mathfrak a_{\vec\ell,\mu}^{\delta_{0}}]f\|_p^p)^{\frac 1p}.
\]
Now, Proposition \ref{prop:avl} gives that
\[
(\sum_{z\in \mathcal I_\vl}\|\mathcal A[\mathfrak a_{\vec\ell,z}^{\delta_{\vec\ell}}]f\|_p^p)^{\frac 1p}\lesssim \lambda^{d\epsilon_1}(\frac{\delta_{ 0}}{\delta_{1}})^{d-1-\frac{d+2}{p}+\epsilon_1}(\sum_{\mu\in \delta_{1}\Z^{d-1}}\|\mathcal A[\mathfrak a_{\vec\ell,\mu}^{\delta_{1}}]f\|_p^p)^{\frac 1p}.
\]
Applying Proposition \ref{prop:avl} repeatedly, we have
\[
\begin{aligned}
(\sum_{z\in \mathcal I_\vl}\|\mathcal A[\mathfrak a_{\vec\ell,z}^{\delta_{\vec\ell}}]f\|_p^p)^{\frac 1p}&\lesssim
\lambda^{d\epsilon_1}(\frac{\delta_{0}}{\delta_{N_0}})^{d-1-\frac{d+2}p+\epsilon_1}(\sum_{\tilde z\in \delta_{N_0}\Z^{d-1}}\|\mathcal A[\mathfrak a_{\vec\ell,\tilde z}^{\delta_{N_0}}]f\|_p^p)^{\frac 1p},
\end{aligned}
\]
 since $N_0$ is bounded above by a constant.
Considering  $(\delta_0/\delta_{N_0})=\lambda^{-2\epsilon_1}\delta_\vl2^{\ell_2}\le \lambda^{1/2}$, we get the desired.
\end{proof}

In the rest of this section, we verify Proposition \ref{prop:avl}.

\subsection{Fourier support of $\mathcal A[\mathfrak a_{\vl,\tilde \mu}^{\delta_{n+1}}]$}
In this subsection, we set $2^{\ell_3}=\lambda^{1/2}$ when $k=2$.

For proving Proposition \ref{prop:avl}, we analyze the behavior of the Fourier support of $\mathcal A[\mathfrak a_{\vl,\tilde \mu}^{\delta_{n+1}}]$ for $\tilde \mu\in\delta_{n+1}\Z^{d-1}\cap B(\mu,2\delta_n)$. As a preliminary step, we first check key properties of $\gamma$.

Fix  $0\le n\le N_0-1$ and $\mu\in \delta_{n}\Z^{d-1}$ satisfying $S_\vl\cap B(\mu,\delta_n)\neq\emptyset$.
By Proposition \ref{prop:Dmg},  we can obtain
\Be
\label{ineq:localgamma}
\mathfrak D^2\gamma(u)\sim \lambda^{-\frac 12}2^{\ell_2} \text{ and }\mathfrak D^3\gamma(u)\lesssim \lambda^{-\frac 12}2^{\ell_3},
\Ee
for $u\in B(\mu,2\delta_n)$, if $\lambda$ is sufficiently large.
Taylor's theorem tells that
\Be
\nonumber
\gamma(u)=\gamma(\mu)+\nabla\gamma(\mu)\cdot(u-\mu)+
\frac 12\langle H\gamma(\mu)(u-\mu),u-\mu\rangle+ O(\lambda^{-\frac 12}2^{\ell_3}\delta_{n}^3),
\Ee
where $H\gamma$ is the Hessian matrix of $\gamma$.
By definition of $\delta_n$, one can check that
\Be
\label{ineq:deln}
\lambda^{- 1/2}2^{\ell_3}\delta_{n}^3\lesssim \lambda^{-1/2} 2^{\ell_2}\delta_{n+1}^2.
\Ee
This suggests that the term $O(\lambda^{-\frac 12}2^{\ell_3}\delta_{n}^3)$ is sufficiently small. 
The properties of $\gamma$ are determined by the main term, $\langle H\gamma(\mu)(u-\mu),u-\mu\rangle$.
Let $\{e^\mu_j\}_{j=1}^{d-1}$ be orthonormal eigenvectors of $H\gamma(\mu)$, where the eigenvalue corresponding to $e_1^\mu$ has the largest absolute value among all eigenvalues of $H\gamma(\mu)$. 
By \eqref{ineq:localgamma}, it holds that
\Be
\label{ineq:e1curvature}
|e_1^\mu\cdot H\gamma(\mu)e_1^\mu|\sim \lambda^{-1/2} 2^{\ell_2}.
\Ee
This indicates that the hypersurface $\Gamma$ exhibits nonvanishing principal curvature in the direction $e_1^\mu$. However, for the remaining directions in $\{e_1^\mu\}^\perp$, no specific curvature conditions can be ensured for $\Gamma$. 
Therefore, we focus on the behavior of $\gamma$ along the direction $e_1^\mu$.

Define 
\[
\mathcal J_\mu(\kappa):=\{\tilde \mu\in \delta_{n+1}\Z^{d-1}\cap B(\mu,2\delta_n):|\Pi_{e_1^\mu}( \tilde \mu)-\kappa|\le\delta_{n+1}\},
\]
for $\kappa \in \{e_1^\mu\}^\perp$. We can observe that
\Be
\label{ineq:Jmukappa}
\{\tilde \mu\in \delta_{n+1}\Z^{d-1}:\tilde \mu\in B(\mu,2\delta_{n})\}\subset \cup_{\kappa\in \mathcal K_\mu}\mathcal J_\mu(\kappa) 
\Ee
where $\mathcal K_\mu:=\{\kappa=\sum_{j=2}^{d-1}a_je_j^\mu: a_j\in \delta_{n+1}\Z \text{ and } \kappa\in \Pi_{e_1^\mu}(B(\mu,2\delta_n)) \}$. Note that
\Be
\label{ineq:sizeK}
| \mathcal K_\mu|\lesssim (\delta_{n}/\delta_{n+1})^{d-2}.
\Ee
For $\kappa\in \mathcal K_\mu$, choose an element $\tilde \mu_\kappa\in \mathcal J_\mu(\kappa)$. And set
\[
\gamma_{\kappa}(u):=\gamma_{\kappa}(\ell_2,\lambda)(u)=\lambda^{1/2}2^{-\ell_2}(\gamma(u)-\gamma(\tilde \mu_\kappa)-\nabla\gamma(\tilde \mu_\kappa)\cdot(u-\tilde \mu_\kappa)),
\]
Then we can check the following properties of $\gamma_\kappa$.
\begin{prop}\label{prop:gammakappa}
    Let $\kappa\in \mathcal K_\mu$ and  $u\in \cup_{\tilde\mu\in \mathcal J_\mu(\kappa)}B(\tilde \mu,\delta_{n+1})$. Then there exists a constant $c_\kappa$ such that $c_\kappa\sim 1$,
  \[
        \gamma_{\kappa}(u)=\frac {c_\kappa}2((u-\tilde \mu_\kappa)\cdot e_1^\mu)^2+ O(\delta_{n+1}^2),
        \]
        and  \[
        \nabla\gamma_{\kappa}(u)=c_\kappa(u-\tilde \mu_\kappa)\cdot e_1^\mu e_1^\mu+ O(\delta_{n+1}).
        \]
\end{prop}

\begin{proof}
We prove the first part of the proposition. The second part follows from almost the same argument, so we omit the details.

First, we can observe that for $u\in \cup_{\tilde\mu\in \mathcal J_\mu(\kappa)}B(\tilde \mu,\delta_{n+1})$,
\begin{equation}
\label{u-num}
    |u-\tilde\mu_\kappa|\lesssim \delta_{n} \text{ and }|\Pi_{e_1^\mu}(u-\tilde\mu_\kappa)|\lesssim \delta_{n+1}.
\end{equation}
By Taylor's theorem combined with \eqref{u-num} and \eqref{ineq:deln}, we can obtain that 
\[
\gamma_\kappa(u)=\frac {\lambda^{1/2}2^{-\ell_2}}2\langle H\gamma(\tilde\mu_\kappa)(u-\tilde\mu_\kappa),u-\tilde\mu_\kappa\rangle+ O(\delta_{n+1}^2).
\]
Define $c_\kappa=\lambda^{1/2}2^{-\ell_2}\langle H\gamma(\tilde\mu_\kappa)e_1^\mu,e_1^\mu\rangle$. By \eqref{ineq:e1curvature}, we can see that $c_\kappa\sim 1$.
Since $u-\tilde\mu_\kappa=\sum_j (u-\tilde\mu_\kappa)\cdot e_j^\mu e_j^\mu$, it is sufficient to show that
\[
\big|\sum_{j=2}^{d-1}\langle H\gamma(\tilde\mu_\kappa)e_j^\mu,e_j^\mu\rangle((u-\tilde\mu_\kappa)\cdot e_j^\mu)^2\big|\lesssim \lambda^{-1/2}2^{\ell_2}\delta_{n+1}^2.
\]
It follows from \eqref{ineq:localgamma} and \eqref{u-num}. 
\end{proof}

Using Proposition \ref{prop:gammakappa}, we show that after applying a suitable linear transform, the Fourier supports of $\mathcal A[\mathfrak a_{\vl,\tilde \mu}^{\delta_{n+1}}]f$ for $\tilde \mu\in \mathcal J_\mu(\kappa)$ are contained in a neighborhood of a cone.
For this purpose, define a linear transformation $T_{\kappa}=T_\kappa(\ell_2,\lambda):\R^{d+1}\rightarrow \R^{d+1}$ as follow:
\[
T_{\kappa}(\tau,\xi)=\lambda^{-1/2}2^{-\ell_2}(\tau+\xi\cdot(\tilde\mu_\kappa,\gamma(\tilde\mu_\kappa)),\xi'+\xi_{d}\nabla\gamma(\tilde\mu_\kappa),\lambda^{-1/2}2^{\ell_2}\xi_{d})
\]
where $\xi=(\xi',\xi_d)$. And define
\[
\Gamma_{\kappa}(u)=(u-\tilde\mu_\kappa,\gamma_\kappa(u)).
\]
Then by definition, we can observe that
\Be
\label{eq:Tkappa}
(\tau,\xi)\cdot(1,\Gamma(u))=\lambda^{\frac 12}2^{\ell_2}T_{\kappa}(\tau,\xi)\cdot(1,\Gamma_{\kappa}(u)).
\Ee
Utilizing Proposition \ref{prop:gammakappa}, we can show that the projection of $T_{\kappa}(\supp \widehat{\mathcal A[\mathfrak a_{\vl,\tilde \mu}^{\delta_{n+1}}]f})$ onto the space generated by $\{(1,\vec 0,0),(0,e_1^\mu,0),(0,\vec 0,1)\}$ is in a neighborhood of a cone. 
Indeed, we have the following:
\begin{prop}\label{prop:projcone}
Let $s_{\tilde \mu}=e_1^\mu\cdot(\tilde \mu-\tilde \mu_\kappa)$ for $\tilde \mu\in \mathcal J_\mu(\kappa)$. Suppose that $(\tau,\xi)$ is contained in the Fourier support of $\mathcal A[\mathfrak a_{\vl,\tilde \mu}^{\delta_{n+1}}]f$.  Then for $m=0,1,2$
    \Be
    \label{projcone}
| T_{\kappa}(\tau,\xi)\cdot\tilde{\mathbf r}_{\kappa}^{(m)}(s_{\tilde \mu})|\lesssim \delta_{n+1}^{2-m},
\Ee
and
\[
| T_{\kappa}(\tau,\xi)\cdot\tilde{\mathbf r}_{\kappa}^{(2)}(s_{\tilde \mu})|\sim 1,
\]
where
\[
\tilde{\mathbf r}_{\kappa}(s):=(1,se_1^\mu,\frac{c_\kappa s^2}{2}).
\]
\end{prop}
\begin{proof}
When $m=2$, one can easily check that $ T_{\kappa}(\tau,\xi)\cdot\tilde{\mathbf r}_{\kappa}^{(2)}(s_{\tilde \mu})=c_\kappa\lambda^{-1}\xi_d$. Thus, the fact $\xi_d\sim \lambda$, which can be deduced from \eqref{maincaes}, implies the desired result for $m=2$. 

For the cases $m=0,1$, we first observe that \eqref{maincaes} holds for some $u\in B(\tilde \mu,\delta_{n+1})$, due to the support condition of $\mathfrak a_{\vl,\tilde\mu}^{\delta_{n+1}}$.
By Taylor's theorem, the condition \eqref{maincaes} for $u\in B(\tilde \mu,\delta_{n+1})$ implies the following estimates:
\[
|(\tau,\xi)\cdot (1,\Gamma(\tilde \mu))|\lesssim \lambda^{\epsilon_1}+\lambda^{\frac 12+\epsilon_1}\delta_{n+1}+\lambda^{\frac 12}2^{\ell_2}\delta_{n+1}^2,
\]
\Be
\nonumber
|\xi\cdot\partial_j\Gamma(\tilde\mu)|\lesssim \lambda^{\frac 12+\epsilon_1}+\lambda^{\frac 12}2^{\ell_2}\delta_{n+1},
\Ee
for $j=1,\cdots,d-1$.
Then the relation \eqref{eq:Tkappa} yields that
\Be
\label{ineq:Tkg}
|T_{\kappa}(\tau,\xi)\cdot\partial^\alpha(1,\Gamma_{\kappa}(\tilde \mu))|\lesssim \delta_{n+1}^{2-|\alpha|},
\Ee
for $|\alpha|\le 1$. Here, we use the condition
\Be
\label{cond:deltan}
\lambda^{\epsilon_1}2^{-\ell_2}\le \delta_{n+1},
\Ee
together with the fact that $\lambda^{-1/2}\lesssim \delta_{n+1}$, which follows from \eqref{cond:lm}.

Now, consider the case $m=0$. By Proposition \ref{prop:gammakappa}, we can  observe that 
\[
(1,\Gamma_\kappa(\tilde \mu))=\tilde{\mathbf r}_\kappa(s_{\tilde \mu})+\sum_{j=2}^{d-1}(\tilde \mu-\tilde\mu_\kappa)\cdot e_j^\mu (0,e_j^\mu,0) +O(\delta_{n+1}^2).
\]
Considering \eqref{ineq:Tkg}, $| T_{\kappa}(\tau,\xi)\cdot\tilde{\mathbf r}_{\kappa}(s_{\tilde \mu})|$ is bounded by
\[
| \sum_{j=2}^{d-1}(\tilde \mu-\tilde\mu_\kappa)\cdot e_j^\mu T_{\kappa}(\tau,\xi)\cdot(0,e_j^\mu,0)|+O(\delta_{n+1}^2).
\]
By \eqref{u-num} which explains that $|(\tilde \mu-\tilde \mu_\kappa)\cdot e_j^\mu|\lesssim \delta_{n+1}$ for $j=2,\cdots,d-1$, 
verifying  \eqref{projcone} for $m=0$ reduces to showing that, for $j=2,\cdots,d-1$,
\Be
\label{ineq:Tksupport}
|T_{\kappa}(\tau,\xi)\cdot(0,e_j^\mu,0)|\lesssim \delta_{n+1}.
\Ee
One can observe that $T_{\kappa}(\tau,\xi)\cdot(0,e_j^\mu,0)$ is equal to 
\[
T_{\kappa}(\tau,\xi)\cdot(0,e_j^\mu,e_j^\mu\cdot \nabla\gamma_\kappa(\tilde \mu))-\lambda^{-1}\xi_de_j^\mu\cdot \nabla\gamma_\kappa(\tilde \mu).
\]
By \eqref{ineq:Tkg} and the definition of $\gamma_\kappa$, to prove \eqref{ineq:Tksupport}, it is sufficient to show that 
 $ |e_j^\mu\cdot (\nabla\gamma(\tilde \mu)-\nabla\gamma(\tilde \mu_\kappa))|\lesssim \lambda^{-1/2}2^{\ell_2}\delta_{n+1}$.
Using Taylor's theorem, we can write 
\[
\nabla\gamma(\tilde \mu)-\nabla\gamma(\tilde \mu_\kappa)=H\gamma(\mu)(\tilde \mu-\tilde \mu_\kappa)+O(\lambda^{-\frac 12}2^{\ell_3}\delta_{n}^2).
\]
Since $e_j^\mu$ is the eigenvector of $H\gamma(\mu)$,  the desired inequality follows from  \eqref{u-num}, which tells that $|(\tilde\mu-\tilde \mu_\kappa)\cdot e_j^\mu|\lesssim \delta_{n+1}$, combined with \eqref{ineq:deln}, which implies that $\lambda^{-1/2}2^{\ell_3}\delta_n^2\lesssim \lambda^{-1/2}2^{\ell_2}\delta_{n+1}$.

When $m=1$, Proposition \ref{prop:gammakappa} says that
\[
(0,(e_1^\mu\cdot \nabla) \Gamma_\kappa(\tilde\mu))= \tilde{\mathbf r}_\kappa'(s_{\tilde \mu}) +O(\delta_{n+1}).
\]
Thus, \eqref{projcone} for $m=1$ follows from \eqref{ineq:Tkg}.
\end{proof}
Now, we are ready to prove Proposition \ref{prop:avl}. 
\subsection{Proof of Proposition \ref{prop:avl}}
Recall \eqref{ineq:Jmukappa}. Since $\mathcal J_\mu(\kappa)$ is finite set, there exist collections $\widetilde{\mathcal J}_\mu(\kappa)\subset \mathcal J_{\mu}(\kappa)$ such that $\widetilde{\mathcal J}_\mu(\kappa)\cap \widetilde{\mathcal J}_\mu(\kappa')=\emptyset$ for $\kappa,\kappa'\in \mathcal K_\mu$  satisfying $\kappa\neq \kappa'$ and
\[
\{\tilde \mu\in \delta_{n+1}\Z^{d-1}:\tilde \mu\in B(\mu,2\delta_{n})\}\subset \cup_{\kappa\in \mathcal K_\mu}\widetilde{\mathcal J}_\mu(\kappa). 
\]
Then
by H\"older's inequality, one can find that
\[
\begin{aligned}
    \|\mathcal A[\mathfrak a_{\vec\ell,\mu}^{\delta_{n}}]f\|_p&\lesssim 
\|\sum_{\kappa\in\mathcal K_\mu
 }\sum_{\tilde \mu\in \widetilde{\mathcal J}_\mu(\kappa)}\mathcal A[\mathfrak a_{\vec\ell,\tilde z}^{\delta_{n+1}}]f\|_p \\
 &\lesssim (\frac{\delta_n}{\delta_{n+1}})^{(d-2)(1-\frac1p)}(\sum_{\kappa\in\mathcal K_\mu
 }\|\sum_{\tilde \mu\in \widetilde{\mathcal J}_\mu(\kappa)}\mathcal A[\mathfrak a_{\vec\ell,\tilde \mu}^{\delta_{n+1}}]f\|_p^p)^{\frac 1p}.
\end{aligned}
\]
Here, for the second line, we use \eqref{ineq:sizeK}. Since 
\[
(\sum_{\mu\in \delta_n\Z^{d-1}}\sum_{\kappa\in\mathcal K_\mu
 }\sum_{\tilde \mu\in \widetilde{\mathcal J}_\mu(\kappa)}\|\mathcal A[\mathfrak a_{\vec\ell,\tilde \mu}^{\delta_{n+1}}]f\|_p^p)^{\frac 1p}\lesssim (\sum_{\tilde \mu \in \delta_{n+1}\Z^{d-1}
 }\|\mathcal A[\mathfrak a_{\vec\ell,\tilde \mu}^{\delta_{n+1}}]f\|_p^p)^{\frac 1p},
\]
Proposition \ref{prop:avl} follows from
\[
\|\sum_{\tilde \mu\in\widetilde{\mathcal J}_\mu(\kappa)}\mathcal A[\mathfrak a_{\vec\ell,\tilde z}^{\delta_{n+1}}]f\|_p\lesssim (\frac{\delta_{n}}{\delta_{n+1}})^{1-\frac 4p+\epsilon}(\sum_{\tilde \mu\in \widetilde{\mathcal J}_\mu(\kappa)}\|\mathcal A[\mathfrak a_{\vec\ell,\tilde \mu}^{\delta_{n+1}}]f\|_p^p)^{\frac 1p}.
\]    
It is a direct consequence of Theorem \ref{BD decoup} combined with Proposition \ref{prop:projcone} and scaling.
\qed

\subsection*{Acknowledgments} This research was supported by a KIAS Individual Grant SP089101.

\end{document}